\documentclass[12pt,a4paper,reqno]{article}
\usepackage[scaled]{helvet}

\usepackage{graphicx}
\usepackage{epstopdf} 
\usepackage{amsfonts}
\usepackage{amsthm}
\usepackage{amsmath}
\usepackage{amscd}
\usepackage[mathscr]{eucal}
\usepackage{indentfirst}

\usepackage{pict2e}
\usepackage{epic}
\numberwithin{equation}{section}
\usepackage{MnSymbol}
\usepackage[top=1.5cm, bottom=1.5cm, left=1.5 cm, right=1cm]{geometry}
\usepackage{tikz}
\usepackage{enumitem}

\usepackage{hyperref}
\numberwithin{equation}{section}

\theoremstyle{plain}
\newtheorem{Th}{Theorem}[section]
\newtheorem{Lemma}[Th]{Lemma}
\newtheorem{Cor}[Th]{Corollary}

 \theoremstyle{definition}
\newtheorem{Def}[Th]{Definition}

\newtheorem{Rem}[Th]{Remark}
\newtheorem{?}[Th]{Problem}

\newcommand{\R}{\mathbb{R}}

\newcommand{\F}{{\mathcal F}^{\mu}}

\begin{document}
\vskip .2cm 
\begin{center}
\huge\textbf{Donoho-Stark's Uncertainty Principles in Real Clifford Algebras}
\end{center}

\begin{center}
\textbf{Youssef El Haoui$^{1,}$\footnote[1]{Corresponding author.}, Said Fahlaoui$^1$}

 $^1$Department of Mathematics and Computer Sciences, Faculty of Sciences, University Moulay Ismail, Meknes 11201, Morocco\\

E-MAIL: y.elhaoui@edu.umi.ac.ma, s.fahlaoui@fs.umi.ac.ma
\end{center}

\vspace{1cm}
\begin{abstract}
The Clifford Fourier transform (CFT) has been shown to be a powerful tool in the Clifford analysis. In this work, several uncertainty inequalities are established in the real Clifford algebra $Cl_{(p,q)}$, \ including the Hausdorf-Young inequality, and three qualitative uncertainty principles of Donoho-Stark. 
\end{abstract}

\textbf{Key words:} Clifford algebras, Clifford-Fourier transform, Uncertainty principle, Donoho-Stark's uncertainty principle.


\section{Introduction}
 It is well known that the uncertainty principles (UPs) give information about a function and its Fourier transform. Their importance is due to their applications in different areas, e.g. quantum physics and signal processing. In quantum physics, they tell us that the position and the momentum of a particle cannot both be measured with precision.

 The qualitative UP is a kind of UPs, which tells us how a  signal $f$\ and its Fourier transform $\hat{f}$,\ behave under certain conditions. One such example can be the Donoho-Stark's UP \cite{DS89}, which expresses the limitations on the simultaneous concentration of $f$, and $\hat{f}$.

The aim of this work is to generalize Donoho-Stark's UP in Clifford's analysis, using the basic properties of Clifford's algebras and its Fourier transform.

For more details on Clifford Fourier's transformations, their historical development and applications, we refer to \cite{BDS05, BDS09, HM08, HI12}.

In \cite{CKL15} Thm. 5.1, and \cite{LI18} Thm. 8, the authors establish , in different ways, the UP of Donoho-Stark in quaternion algebra which is isomorphic at  ${Cl}_{(0,2)}.$

The first inequality we deal with is a generalization of the Hausdorf-Young inequality by means of the kernel of the CFT introduced by \cite{HI12}.

Based on this inequality, and following  the Donoho-Stark'UP proof techniques for the Dunkel-trabsform \cite{SO14}, we investigate three inequalities in terms of ''$\epsilon$-concentration'' in the Clifford algebra ${Cl}_{(p,q)}.$

This paper is organized as follows. Section 2 is devoted to a reminder of the basics of Clifford algebras. In section 3, we introduce the CFT and review its important properties, and prove the Hausdorf-Young inequality. In section 4, we define the concept of ''$\epsilon$-concentration'' in CFT-domain, and establish UPs of concentration type, then prove Donoho-Stark bandlimited UP for the CFT. Finally, we give a conclusion in  section 5.

%

\section{Preliminaries}

Let $\{e_1,e_2,\dots ,e_n\}$ be an orthonormal basis of the real Euclidean vector space $\R^{(p,q)}$, with $p+q=n.$

The Clifford geometric algebra (see \cite{MU91})   over $\R^{(p,q)}$ denoted by $Cl_{(p,q)}$, is defined as an associative, non commutative algebra which has the graded 2$^n$-dimensional basis

$\{{1,e}_1,e_2,\dots ,e_n,e_{12},e_{13},\ e_{23},\dots ,e_1e_2\dots e_n\}.$

The multiplication of the basis vectors satisfy the rules

 $e_ke_l+e_le_k=2 \epsilon_k{\delta }_{k,l}$, for $1\ \le k,l\le $n,

With ${\delta }_{k,l}$ is the Kronecker symbol, and $\epsilon_k=+1,\ $ for $k=1,..,p,$\ and $\ \epsilon_k=-1,\ $ for $k=p+1,..,n.$

Every element $\ f\ $ of Clifford algebra $Cl_{(p,q)}$, is called multivector, and can be expressed in the form

\[f\left(x\right)=\sum_A{f_A\left(x\right)e_{A,}}\] 

where $f_A,\ $are real-valued functions, $e_A=e_{{\alpha }_1{\alpha }_2\dots {\alpha }_k}=e_{{\alpha }_1}e_{{\alpha }_2}{\dots e}_{{\alpha }_k}$, and$\ 1\le {\alpha }_1\le {\alpha }_2\le \dots \le {\alpha }_k\le n$, with ${\alpha }_i\in \{1,2,\dots ,n\}$, and $e_{\emptyset }=1$.

Also, a multivector $f\in {Cl}_{(p,q)}\ ,$\ can be written as\\

\[f=\sum^{k=n}_{k=0}{{<f>}_k}={<f>}_0+{<f>}_1+{<f>}_2+\dots +{<f>}_n.\]
  where ${<f>}_k=\sum_{\left|A\right|=k}{{\alpha }_A}e_A$, denote the $k-$vector part of$\ f$.

And the reverse $\tilde{f}$ of is given by 

\[\tilde{f}=\sum^{k=n}_{k=0}{{(-1)}^{\frac{k(k-1)}{2}}}\overline{{<f>}_k},\] 

Where $\overline{f}$ means to change in the basis decomposition of $f$ the sign of every vector of negative square  $\overline{e_A}=\epsilon_{{\alpha }_1}e_{{\alpha }_1}\dots \epsilon_{{\alpha }_k}e_{{\alpha }_k},\ \ 1\le {\alpha }_1\le {\alpha }_2\le \dots \le {\alpha }_k\le n\ $.

For $f,\tilde{g}\in {Cl}_{(p,q)}$, the scalar product $f*\tilde{g}$, is defined by

\[f*\tilde{g}={<f\tilde{g}>}_0 =\sum_A{f_Ag_A\ }.\] 
 
 In particular, if $f=g$, then we obtain the modulus of a multivector $f\in {Cl}_{(p,q)}$, defined as 

\begin{equation} \label{Cliff_norm}
 \left|f\right|=\sqrt{f*\tilde{f}}=\sqrt{\sum_A{f^2_A}}.
\end{equation}

For $1\le a<\infty $ , The linear spaces ${\ L^a(\R^{(p,q)},Cl}_{(p,q)})$ are introduced as :

${ L^a(\R^{\left(p,q\right)},Cl}_{\left(p,q\right)})=\{\R^{\left(p,q\right)}\to {Cl}_{\left(p,q\right)}:\ {\left\|f\right\|}_a={\left(\int_{\R^{\left(p,q\right)}}{{\left|f\left(t\right)\right|}^a}dt\right)}^{\frac{1}{a}}<\infty \}$.

 For $a=\infty ,$ the $L^{\infty }$- norm is defined by 
            
\[{\left\|f\right\|}_{\infty }=ess\ {sup}_{t\in \R^{\left(p,q\right)|}}\left|f\left(t\right)\right|.\] 

\begin{Lemma}\label{Brax}

For $\lambda ,\rho \in{Cl}_{(p,q)}$, the following property hold 

\[\left|\rho \lambda \ \right|\le 2^n\left|\lambda \right|\left|\rho \right|.\] 
\end{Lemma}

\begin{Lemma}\label{Eul}

Let  $\theta \in \R,$\ and $\mu \in {Cl}_{(p,q)}$ , with ${\mu }^2=-1,\ \ $we have a natural generalization of Euler's formula for Clifford algebra, as follows 

\[e^{\theta \mu}={\cos  \left(\theta \right) }+\mu {\sin  \left(\theta \right)\ }.\] 
\end{Lemma}

Proof. As ${\mu}^2=-1,$\ we have for any real $\theta$
\begin{eqnarray*}
e^{\mu \theta} &=& \sum^{\infty }_{k=0}{\frac{{(\mu \theta )}^k}{k!}} \\
&=& \sum^{\infty }_{k=0}{{(-1)}^k\frac{{\theta }^{2k}}{(2k)!}}+\mu \sum^{\infty }_{k=0}{{(-1)}^k\frac{{\theta }^{2k+1}}{(2k+1)!}}\\
&=& \cos \left(\theta \right)+\mu \sin  \left(\theta \right) .\hspace{7cm}\square
\end{eqnarray*}

\section{Clifford-Fourier transform}
In this section, we introduce the Clifford Fourier transform (CFT),recall its properties, add one result related to the kernel of the CFT, and we prove the Hausdorf-Young inequality  associated with the CFT.

\begin{Def}

Let  $\mu \in{\ Cl}_{(p,q)}$ be a square root of -1, i.e. ${\mu }^2=-1.$

The general Clifford Fourier transform (CFT) (see \cite{HI12}) of $f\in {\ L^1(\R^{\left(p,q\right)},Cl}_{\left(p,q\right)})$, with respect to $\mu $\ is

\[\F\left\{f\right\}(\xi ) = \int_{\R^{(p,q)}}{f(t){e^{-\mu u(t,\xi )}dt}}.\] 

Where $dt={dt}_1\dots {dt}_n,\ t,\xi \in \R^{\left(p,q\right)}$, and $u:\ \R^{\left(p,q\right)}\times \R^{\left(p,q\right)}\to \R.$\\
 We assume, in the rest of this work, that 
\[u\left(t,\xi \right)=\sum^{l=n}_{l=1}{t_l{\xi }_l}.\]

\end{Def}
\subsection{Properties of the CFT}
In the following, we give some important properties of the CFT, For more detailed discussions of the properties of the CFT and their proofs, see e.g. \cite{BDS05,HI12,HM08}

\hspace*{0.6 cm} $ \tikz\draw[black,fill=black] (0,0) circle (.5ex);$\ Left Linearity

For $f_1,\ f_2\in{\ L^1(\R^{\left(p,q\right)},Cl}_{\left(p,q\right)})$,\ and constants $\alpha , \beta \in {Cl}_{\left(p,q\right)}$, 

\[ \F\left\{\alpha f_1+\beta f_2\right\}= \alpha \F\left\{f_1\right\}+\ \beta \F\left\{f_2\right\}.\] 

\hspace*{0.6 cm} $ \tikz\draw[black,fill=black] (0,0) circle (.5ex);$\ Inversion formula

For $f,\ \F\left\{f\right\}\in{L^1(\R^{\left(p,q\right)},Cl}_{\left(p,q\right)})$, we have
\begin{equation}\label{invers}
f\left(t\right)=\frac{1}{{(2\pi )}^n}\int_{\R^{(p,q)}}{\F\left\{f\right\}\left(\xi \right){e^{\mu u\left(t,\xi \right)}d\xi .}}
\end{equation}
Where $d\xi ={d\xi }_1$\dots ${d\xi }_n$, $t,\xi \in \R^{\left(p,q\right)}$.

\hspace*{0.6 cm} $ \tikz\draw[black,fill=black] (0,0) circle (.5ex);$\ For the function $f\ \in{\ L^2(\R^{\left(p,q\right)},Cl}_{\left(p,q\right)}$, one has the Parseval identity
\begin{equation}\label{pars}
{\left\|\F\left\{f\right\}\right\|}_2={(2\pi )}^{\frac{n}{2}}{\left\|f\right\|}_2.
\end{equation}

\begin{Lemma}\label{Jday}
For $x,y \in \R^{(p,q)}$, and $\mu \in {Cl}_{(p,q)},\ $ with ${\mu }^2=-1,$ the following inequality holds: 
\[  |e^{-\mu u(x,y)}| \le {(1+{|\mu |}^2)}^{\frac{1}{2}}.\]
\end{Lemma}

Proof. By means of Lemma \ref{Eul} and the definition \eqref{Cliff_norm} of the Clifford norm, we obtain 
\begin{eqnarray*}
{\left|e^{-\mu u\left(t,\xi \right)}\right|}^2 &=& {\left|{\cos\left(u\left(t,\xi \right)\right)}\right|}^2+\sum_A{{\left|{\sin \left(u\left(t,\xi \right)\right)}\right|}^2{\left|{\mu }_A\right|}^2}\\
&\le& 1+\sum_A{{\left|{\mu }_A\right|}^2}=1+{|\mu |}^2.                             
\end{eqnarray*}

Therefore,

\hspace*{4 cm}${\left|e^{-\mu u\left(t,\xi \right)}\right|}\le {({\rm 1+}{|\mu |}^2)}^{\frac{1}{2}}\hfill \square$ 

However  , by combining lemma \ref{Jday} and lemma \ref{Brax}, we do have the following result:

\begin{Lemma}

Let $\lambda \in {Cl}_{(p,q)}$, and $\mu \in{Cl}_{(p,q)}$ be a square root of -1.

Then,
\begin{equation}\label{usf}
|\lambda e^{-\mu u(x,y)}|\le 2^n\left|\lambda \ \right|{(1+{|\mu |}^2)}^{\frac{1}{2}}.
\end{equation}
\end{Lemma}

\begin{Th}{Hausdorf-Young inequality associated with CFT}

Let ${\ f\in L^a(\R^n,Cl}_{(p,q)}), 1\le a\le 2$,

 then  ${\ \F\left\{f\right\}\in L^b(\R^n,Cl}_{(p,q)})\ \ $with $\frac{1}{a}+\frac{1}{b}=1$, 

and we have 
\begin{equation} \label{haus}
{\left\|\F\left\{f\right\}\right\|}_b\le {\left(2^nC_{\mu }\right)}^{1-\frac{2}{b}}{\left(2\pi \right)}^{\frac{2n}{b}}{\left\|f\right\|}_a,
\end{equation}
\ \ where  $C_{\mu }={(1+{|\mu |}^2)}^{\frac{1}{2}}$.
\end{Th}
Proof. 
We have 
\begin{eqnarray*}
\left|\F\left\{f\right\}(\xi )\right|&=&\left|\int_{\R^n}{f(t){e^{-\mu u(t,\xi )}dt}}\right|\\
&\le& \int_{\R^n}{\left|f(t)e^{-\mu u(t,\xi )}\right|dt}\\ 
&\le& 2^n{(1+{|\mu |}^2)}^{\frac{1}{2}}\int_{\R^n}{\left|f(t)\right|dt}.                             
\end{eqnarray*}
Where  we used \eqref{usf}.\\
Thus, \begin{equation}\label{infini_norm}
{\left\|\F\left\{f\right\}\right\|}_{\infty }\le 2^nC_{\mu }{\left\|f\right\|}_1.
\end{equation}
Hence,  the CFT is of type (1,$\infty $) with norm 2$^n{(1+{|\mu |}^2)}^{\frac{1}{2}}$ ,

On the other hand, by (Parseval ) one sees that the CFT is of type (2,2) with norm ${\left(2\pi \right)}^n$.

We obtain consecutively by the Riesz--Thorin theorem (\cite{LP04}, Thm. 2.1), that the CFT is also of type ($a,b)$, with norm $M_{\theta },$\ such that$\ M_{\theta }\le {(2^nC_{\mu })}^{1-\theta }{({\left(2\pi \right)}^n)}^{\theta }.$

With\ $\frac{1}{a}=\frac{1-\theta }{1}+\frac{\theta }{2}=1-\frac{\theta }{2}$  and  $\frac{1}{b}=\frac{1-\theta }{\infty }+\frac{\theta }{2}=\frac{\theta }{2},$ with 0$\le \theta \le $1.

Then $\frac{1}{a}=1-\frac{1}{b}$, and   1$\le a\le 2$,and \ $M_{\theta }\le {\left(2^nC_{\mu }\right)}^{1-\frac{2}{b}}{\left(2\pi \right)}^{\frac{2n}{b}}.$

This completes the proof.

\section{Donoho-Stark Uncertainty Principles in Clifford algebra $Cl_{(p,q)}$}
The Donoho-Stark UP relies on the concept of $\epsilon-$concentration. We start by giving this definition.\\
Let $T,\ \Omega $ be  measurable subsets of $\R^n$.
And denote by ${P}_{T}, and\ \ Q_{\Omega }$ respectively the time limiting operator, and the Dunkel integral operator given by

\[{\ \ \ \ \ \ \ \ \ \ \ \ \ \ \ \ \ \ P}_{T} f={\chi }_{T}.f, \ \ \ \ \ \ \ \ \ \ \ \ \ \ \ \ \ \ \ \ \ \ \ \ \ \ \ \ \ \F\left(Q_{\Omega }f\right)={\chi }_{\Omega }(\F(f)).\] 
\begin{Def}
 A function $f$ is ${\varepsilon }_T$ concentrated on $T$,\  in $L^a-$norm,

If \\
\[{(\int_{\R^n\backslash T}{{|f(t)|}^adt})}^{\frac{1}{a}}\le {\varepsilon }_T{\left\|f\right\|}_a.\]
\end{Def}

\begin{Rem}\label{Rem1}
If  ${\varepsilon }_T=0,$ then $T$ is the exact support of $f$.
\end{Rem}

\begin{Lemma}

If $1<a\le 2,\ \frac{1}{a}+\frac{1}{b}=1$, and ${\ f\in L^a(\R^n,Cl}_{(p,q)}),$\\
 Then
\begin{equation}\label{up10}
\F\left(Q_{\Omega }P_T\right)\le {2^nC}_{\mu }{|T|}^{\frac{1}{b}}{|\Omega |}^{\frac{1}{b}}{\left\|f\right\|}_a.\end{equation}

Where $|T|$\ is denoted as the Lebesgue measure of $T$.
\end{Lemma}
Proof. Without loss of generality, we assume that $|T|<\infty \ and\ \left|\Omega \right|<\infty \ .$

We have\\
 \[\F\left(Q_{\Omega }P_T\right)={\chi }_{\Omega }\left(\F\left(P_Tf\right)\right),\]
Thus,\\
 \begin{equation}\label{up11}
 {\left\|\F\left(Q_{\Omega }P_T\right)\right\|}_b ={\left(\int_{\Omega }{{\left|\F\left(P_Tf\right)\left(\xi \right)\right|}^b}d\xi \right)}^{\frac{1}{b}}.
\end{equation}
In view of 
\[\F\left(P_Tf\right)\left(\xi \right)=\int_T{f(t){e^{-\mu u(t,\xi )}dt}},\]
We have by H\"older inequality, and \eqref{usf})
\begin{eqnarray*}
\F\left(Q_{\Omega}P_T\right)\left(\xi\right)&\le& 2^n{\left(\int_T{{\left|f\left(t\right)\right|}^a}dt\right)}^{\frac1a}{\left(\int_T{{\left|{e}^{-\mu{\mathbf u}(t,\xi)}\right|}^{b}}dt\right)}^{\frac1{b}}\\
&\le& {{\left\|f\right\|}_a2^nC}_\mu{\left|T\right|}^{\frac1{b}}. 
\end{eqnarray*}
Consequently, \eqref{up11} yields 

\[{\left\|\F\left(Q_{\Omega}P_T\right)\right\|}_{b}\le {2^nC}_\mu{\left|T\right|}^{\frac1{b}}{\left|\Omega\right|}^{\frac1{b}}{\left\|f\right\|}_a.\]

We are now in the position to establish the first UP of concentration type.

\begin{Th}
If a non-zero function ${\ f\in L^1\cap L^a(\R^n,Cl}_{(p,q)}),\ $ is $\varepsilon_{T}$ concentrated on  $T{,\ in\ L}^a-$norm, and $\F$\ is\ ${\varepsilon}_{\Omega}$-concentrated on  $\Omega{,\ in\ L}^{b}-$norm, $\frac1a+\frac1{b}=1$.

Then,
\begin{equation}\label{up1}
{\left\|\F\{f\}\right\|}_{b}\le \frac{C_\mu2^n{|T|}^{\frac1{b}}{|\Omega|}^{\frac1{b}}+C_{h}{\ \varepsilon}_{T}}{1-{\ \varepsilon}_{\Omega}}{\left\|f\right\|}_a.
\end{equation}
With$\ C_{h}={\left(2^nC_{\mu}\right)}^{1-\frac2{b}}{\left(2\pi\right)}^{\frac{2n}{b}}$, the constant of Hausdorf-Young inequality.
\end{Th}
Proof. Withouot loss of generality, we may assume that T and $\Omega$ have finite measure.

Then we have 

\[{\left\|f-P_{T}f\right\|}_a={\left(\int_{\R\backslash T}{{\left|f\left(t\right)\right|}^a}\right)}^{\frac1a}\le \varepsilon_{T}{\left\|f\right\|}_a.\] 

Since $\F\{f\}$ is ${i\ \varepsilon}_{\Omega}$-concentrated on  $\Omega,$\ in\ $L^{b}-$norm, we obtain that
\begin{equation}\label{up13}
{\left\|\F\left\{f\right\}-\F\left\{Q_{\Omega}f\right\}\right\|}_{b}\le {\ \varepsilon}_{\Omega}{\left\|\F\left\{f\right\}\right\|}_{b}.
\end{equation}
On the other hand,
\begin{eqnarray*}
{\left\|\F\left\{Q_{\Omega}f\right\}-\F\left\{Q_{\Omega}P_{T}f\right\}\right\|}_{b}&=&{\left\|{{\chi}_{\Omega}}.\F{\mu}\left\{f\right\}-{{\chi}_{\Omega}}.\F{\mu}\left\{P_{T}f\right\}\right\|}_{b}\\
&\le& {\left\|\F\left\{f-P_{T}f\right\}\right\|}_{b}\\
&\le& C_{h}{\left\|f-P_{T}f\right\|}_a.
\end{eqnarray*}

Where we used the the linearity of $\F$\ in the second inequality, and \eqref{haus} in the last.\\
And consequently, by the triangle inequality
\begin{eqnarray*}
{\left\|\F\left\{f\right\}-\F\left\{Q_{\Omega}P_{T}f\right\}\right\|}_{b}&\le& {\left\|\F\left\{f\right\}-\F\left\{Q_{\Omega}f\right\}\right\|}_{b}+{\left\|\F\left\{Q_{\Omega}f\right\}-\F\left\{Q_{\Omega}P_{T}f\right\}\right\|}_{b}\\
&\le& {\ \varepsilon}_{\Omega}{\left\|\F\left\{f\right\}\right\|}_{b}+C_{h}{\left\|f-P_{T}f\right\|}_a 
\end{eqnarray*}
\begin{equation}\label{up12}
\le {\ \varepsilon}_{\Omega}{\left\|\F\left\{f\right\}\right\|}_{b}+C_{h}\varepsilon_{T}{\left\|f\right\|}_a.
\end{equation}
Moreover, again using the triangle inequality, \eqref{up10}, and \eqref{up12}, implies that
\begin{eqnarray*}
{\left\|\F\left\{f\right\}\right\|}_{b}&\le& {\left\|\F\left\{f\right\}-\F\left\{Q_{\Omega}P_{T}f\right\}\right\|}_{b}+{\left\|\F\left\{Q_{\Omega}P_{T}f\right\}\right\|}_{b}\\
&\le&{\ \varepsilon}_{\Omega}{\left\|\F\left\{f\right\}\right\|}_{b}+C_{h}\varepsilon_{T}{\left\|f\right\|}_a+2^nC_{\mu}{\left|T\right|}^{\frac1{b}}{\left|\Omega\right|}^{\frac1{b}}{\left\|f\right\|}_a.
\end{eqnarray*}
Hence,\\
                          \[\left(1-{\ \varepsilon}_{\Omega}\right){\left\|\F\left\{f\right\}\right\|}_{b}\le (C_{\mu}2^n{|T|}^{\frac1{b}}{|\Omega|}^{\frac1{b}}+C_{h}{\ \varepsilon}_{T}){\left\|f\right\|}_a.\] $\hfill \square$

In view of the Parseval identity \eqref{pars}, and \eqref{up1}

\begin{Cor}\label{cor1}
Suppose that ${\ f\in L^2(\R^n,Cl}_{(p,q)}),$ with $f\ne 0,\ $is $\varepsilon_{T}$ concentrated on  $T{,\ in\ L}^2-$norm, and $\F$\ is\ $ {\varepsilon}_{\Omega}$-concentrated on  $\Omega,\ in\ L^2-$norm.

Then, one has

\[{(2\pi)}^{\frac{n}{2}}(1-{\ \varepsilon}_{\Omega})-C_{h}{\varepsilon}_{T}\le C_{\mu}2^n{|T|}^{\frac12}{|\Omega|}^{\frac12}.\] 
\end{Cor}
Choose  ${\ \varepsilon}_{\Omega}={\ \varepsilon}_{T}=0$,in Corollary \ref{cor1}, and use remark \ref{Rem1},

We do have the following result

\begin{Cor}

Suppose that${\ f\in L^2(\R^n,Cl}_{(p,q)}),$ with $ Supp\ f\subseteq T,\ $ and $Supp \ \F \subseteq \Omega.$

Then,
\[\frac{({\frac{\pi}{2}})^n}{{1+|\mu|}^2}\le |T||\Omega|.\]
\end{Cor}
The second concentration UP of Donoho-Stark associated with CFT, is given in the following theorem.

\begin{Th}
Suppose that a non zero function ${\ f\in L^1\cap L^a(\R^n,Cl}_{(p,q)}),\ 1<a\le 2,$ is $\varepsilon_{T}$ concentrated on  $T{,\ in\ L}^1-$norm, and $\F $is\ ${\varepsilon}_{\Omega}$-concentrated on  $\Omega,\ in\ L^{b}-$norm, $\frac1a+\frac1{b}=1$.

Then,

\[{\left\|\F\{f\}\right\|}_{b}\le \frac{{2^nC}_{\mu}{|T|}^{\frac1{b}}{|\Omega|}^{\frac1{b}}}{\left(1-{\ \varepsilon}_{\Omega}\right)(1-{\ \varepsilon}_{T})}{\left\|f\right\|}_a.\] 
\end{Th}
Proof. We assume that $T$and $\Omega$ have finite measure, we have by triangle inequality and \eqref{up13}
\begin{eqnarray*}
{\left\|\F\left\{f\right\}\right\|}_{b}&\le& {\left\|\F\left\{f\right\}-\F\left\{Q_{\Omega}f\right\}\right\|}_{b}+{\left\|\F\left\{Q_{\Omega}f\right\}\right\|}_{b}\\
&\le& {{\ \varepsilon}_{\Omega}\left\|\F\left\{f\right\}\right\|}_{b}+{\left(\int_{\Omega}{{\left|\F\left\{f\right\}\left(\xi\right)\right|}^a}d\xi\right)}^{\frac1a}.
\end{eqnarray*}
Using
\[{\left(\int_{\Omega}{{\left|\F\left\{f\right\}\left(\xi\right)\right|}^a}d\xi\right)}^{\frac1a}\le {\left\|\F\left\{f\right\}\right\|}_{\infty }{|\Omega|}^{\frac1{b}}.\] 
We indeed obtain by \eqref{infini_norm}
\begin{equation}\label{up21}
\left(1-{\ \varepsilon}_{\Omega}\right){\left\|\F\left\{f\right\}\right\|}_{b}\le {|\Omega|}^{\frac1{b}}2^nC_{\mu}{\left\|f\right\|}_1.
\end{equation}

Furthermore, by assuming that $f$ is $\varepsilon_{T}$ concentrated on  $T,$\ in\ $L^1-$norm, we obtain
\begin{eqnarray*}
{\left\|f\right\|}_1&\le& {\left\|f-P_{T}f\right\|}_1+{\left\|P_{T}f\right\|}_1\\
&\le& {\ \varepsilon}_{T}{\left\|f\right\|}_1+{\left(\int_{T}{{\left|f\left(t\right)\right|}}dt\right)} \\
&\le& {\ \varepsilon}_{T}{\left\|f\right\|}_1+{|T|}^{\frac1{b}}{\left\|f\right\|}_a.
\end{eqnarray*}
Where we used the H\"older inequality.

Thus,
\begin{equation}\label{up22}
\left(1-{\ \varepsilon}_{T}\right){\left\|f\right\|}_1\le {\left|T\right|}^{\frac1{b}}{\left\|f\right\|}_a.
\end{equation}

Combining the results of \eqref{up21} and \eqref{up22} yields the desired result. $\hfill \square$

\subsection{Donoho-Stark bandlimited UP for the CFT}

Let \ ${\mathcal B}^a$($\Omega$), $1\le a\le 2$,\ be the set of functions ${\ g\in L^a(\R^n,Cl}_{(p,q)}),$\ that are bandlimited on $\Omega,$\ i.e. \ $Q_{\Omega}g=g.$

A function $f$ is said $\varepsilon-$bandlimited on $\Omega$\ in\ $ L^a-$norm, if there is $g\in {\mathcal B}^a$($\Omega$), with\\
 \[{\left\|f-g\right\|}_a\le \varepsilon\ {\left\|f\right\|}_a.\]

\begin{Lemma}\label{up31}

Let $g\in {\mathcal B}^a$($\Omega$), $1\le a\le 2,\ $then
\[{\left\|P_{T}g\right\|}_a\le \frac{C_{\mu}C_h}{{\pi}^n}{|T|}^{\frac1a}{|\Omega|}^{\frac1a}{\left\|g\right\|}_a.\] 
\end{Lemma}
Proof. We may assume that $\left|T\right|<\infty \ $and  $\left|\Omega\right|<\infty $.

Using inversion formula \eqref{invers}, and the assumption that $\ g\in {\mathcal B}^a$($\Omega$), we get

\[g\left(t\right)=\frac1{{(2\pi)}^n}\int_{\Omega}{\F\left\{g\right\}(\xi)\ {{e}^{\mu{\mathbf u}(\xi,t)}d\xi}}. \] 

By \eqref{usf} and H\"older inequality, one has
\begin{eqnarray*}
\left|g\left(t\right)\right|&\le& \frac{C_{\mu}}{{\pi}^n}{|\Omega|}^{\frac1a}{\left\|\F\left\{g\right\}\right\|}_{b},\ \ \ \frac1a+\frac1{b}=1.\\
&\le& \frac{C_{\mu}C_{h}}{{\pi}^n}{|\Omega|}^{\frac1a}{\left\|g\right\|}_a. \  \ (By\  \eqref{haus})
\end{eqnarray*}
Thus, we have 
\[{\left\|P_{T}g\right\|}_a={\left(\int_{T}{{\left|g\left(t\right)\right|}^a}dt\right)}^{\frac1a}\le \frac{C_{\mu}C_{h}}{{\pi}^n}{|T|}^{\frac1a}{|\Omega|}^{\frac1a}{\left\|g\right\|}_a.\]   $\hfill \square$

\begin{Th}

 Let ${\ f\in L^1(\R^n,Cl}_{(p,q)}),\ $ be $\varepsilon_{T}$ concentrated on $T,$\ in\ $L^a-$norm and  $\varepsilon$- bandlimited on $\Omega,$\ in\ $L^a-$norm, $1\le a\le 2$

Then,

\[\frac{1-{\ \varepsilon}_{\Omega}-{\ \varepsilon}_{T}}{1+{\ \varepsilon}_{\Omega}}\le \frac{C_{\mu}C_{h}}{{\pi}^n}{|T|}^{\frac1a}{|\Omega|}^{\frac1a}.\] 
\end{Th}
Proof. By definition, there exists $g\in {\mathcal B}^a$($\Omega$), such that ${\left\|f-g\right\|}_a\le {\ \varepsilon}_{\Omega}{\left\|f\right\|}_a$.

This leads to
\begin{eqnarray*}
{\left\|P_{T}f\right\|}_a&\le& {\left\|P_{T}g\right\|}_a+{\left\|P_{T}(f-g)\right\|}_a\\ 
&\le& {\left\|P_{T}g\right\|}_a+{\ \varepsilon}_{\Omega}{\left\|f\right\|}_a. 
\end{eqnarray*}
From Lemma \ref{up31}, and the fact ${\left\|g\right\|}_a\le (1+{\ \varepsilon}_{\Omega})\ {\left\|f\right\|}_a$, we get
\begin{equation}\label{up32}
{\left\|P_{T}f\right\|}_a\le \frac{C_{\mu}C_{h}}{{\pi}^n}{|T|}^{\frac1a}{|\Omega|}^{\frac1a}(1+{ \varepsilon}_{\Omega}){\left\|f\right\|}_a+{\varepsilon}_{\Omega}{\left\|f\right\|}_a=[\frac{C_{\mu}C_{h}}{{\pi}^n}{|T|}^{\frac1a}{|\Omega|}^{\frac1a}(1+{\ \varepsilon}_{\Omega})+ {\varepsilon}_{\Omega}]{\left\|f\right\|}_a. 
\end{equation}
On the other hand, as  $f$ is ${\ \varepsilon}_{T}-$concentrated on $T,$\ in\ $L^a-$norm, we have
\begin{eqnarray*}
{\left\|f\right\|}_a&\le& {\left\|f-P_Tf\right\|}_a+{\left\|P_Tf\right\|}_a\\
&\le& {\ \varepsilon }_{T}{\left\|f\right\|}_a{+\left\|P_Tf\right\|}_a.
\end{eqnarray*}
Then,  
\begin{equation}\label{up33}
{\left\|f\right\|}_a\le \frac{1}{1-{\ \varepsilon }_{T}}{\left\|P_Tf\right\|}_a.
\end{equation}

By combining \eqref{up32} and \eqref{up33}, we conclude the proof.$\hfill \square$

\section{Conclusion}
In this paper, we have proven several uncertainty inequalities for the CFT. The first one is the Hausdorf-Young inequality in the Clifford algebra $Cl_{(p,q)}$,\ which we think will be an important tool in the future to prove other geometric inequalities for the CFT.
The other three inequalities are the generalization of UPs of concentration type, they are $L^a(\R^n,Cl_{(p,q)})$\ versions. Two are dependent on signal f. However, the third is independent of the bandlimited signal $f$.

 

\begin{thebibliography}{9}
\bibitem{BDS05}
F. Brackx, N. De Schepper, and F. Sommen, The Clifford Fourier transform, J. Fourier Anal. Appl., 6(11), pp. 668-681 (2005)
 \bibitem{BDS09}
F. Brackx, N. De Schepper, and F. Sommen, The Fourier transform in Clifford analysis, Advances in Imaging and Electron Physics 156 (2009): pp. 55-201.
 \bibitem{CKL15}
L.P. Chen, K.I. Kou, M.S. Liu, Pitt's inequality and the uncertainty principle associated with
thequaternion Fourier transform, J. Math.Anal.Appl.423(2015), pp. 681–700.
 \bibitem{DS89}
 D.L. Donoho, P.B. Stark, Uncertainty principles and signal recovery, SIAM
J. Appl. Math., 49 (1989), pp. 906-931.
\bibitem{HI12}
E. Hitzer, The Clifford Fourier transform in real Clifford algebras. In E. Hitzer,
K. Tachibana (eds.), ”Session on Geometric Algebra and Applications, IKM
2012”, Special Issue of Clifford Analysis, Clifford Algebras and their Applica-
tions, Vol. 2, No. 3, pp. 223-235, (2013). First published in K. Guerlebeck, T.
Lahmer and F. Werner (eds.), electronic Proc. of 19th International Confer-
ence on the Application of Computer Science and Mathematics in Architecture
and Civil Engineering, IKM 2012, Weimar, Germany, 0406 July 2012. Preprint:
http://vixra.org/abs/1306.0130.

\bibitem{HM08}
E. Hitzer, B. Mawardi, Clifford Fourier Transform on Multivector Fields and
Uncertainty Principles for Dimensions n = 2(mod4) and n = 3(mod4). Adv.
Appl. Cliff. Alg. 18 (2008), pp. 715-736.
\bibitem{LI18}
 P. Lian, Uncertainty principle for the quaternion Fourier transform, J. Math. Anal. Appl. (2018), https://doi.org/10.1016/j.jmaa.2018.08.002
 
 \bibitem{LP04}
 F. Linares, G. Ponce, Introduction to Nonlinear Dispersive Equations, Publica\c c\~oes Matem\'aticas, IMPA, Rio de Janeiro, Brazil, 2004.
\bibitem{MH06}
B. Mawardi, E.M. Hitzer, Clifford Fourier Transformation and Uncertainty
Principle for the Clifford Geometric Algebra Cl(3, 0). Adv. Appl. Cliff. Alg. 16
(2006), pp. 41-61. doi:10.1007/s00006-006-0003-x.
 \bibitem{MU91}
 M. Murray, Clifford algebras and Dirac operators in harmonic analysis, cambridge, university press, First published , 1991.

 \bibitem{SO14}
F. Soltani, Lp-uncertainty principles on Sturm-Liouville hypergroups. Acta Math. Hungar. 2014;142: pp. 433-443.

\end{thebibliography}
\end{document}